\documentclass{icmart}
\usepackage{amsmath,amssymb}
\usepackage{verbatim}

\contact[rudnick@post.tau.ac.il]{Raymond and Beverly Sackler School
of Mathematical Sciences, Tel Aviv University, Tel Aviv 69978,
Israel}

\numberwithin{equation}{section}

\newtheorem{thm}{Theorem}[section]

\newcommand{\Z}{{\mathbb Z}} 

\newcommand{\C}{{\mathbb C}}

\newcommand{\FF}{{\mathbb F}}
\newcommand{\fq}{\mathbb{F}_q}

\newcommand{\tr}{\operatorname{tr}}

\newcommand{\Var}{\operatorname{Var}}

\newcommand{\disc}{\operatorname{disc}}

\newcommand{\Li}{\operatorname{Li}}

\newcommand{\conv}{*}

\newcommand{\Sym}{\operatorname{Sym}}

\newcommand{\EN}{\mathcal N} 
\newcommand{\EM}{\mathcal M} 

\newcommand{\divid}{d} 

\newcommand{\inv}{\theta}

\title[Some problems in analytic number theory for polynomials over a finite field]{Some problems in analytic number theory for polynomials over a finite field}

\author[Zeev Rudnick]
{Zeev Rudnick \thanks{The research leading to these results has
received funding from the European Research Council under the
European Union's Seventh Framework Programme (FP7/2007-2013) / ERC
grant agreement n$^{\text{o}}$ 320755 .}}

\begin{document}

\begin{abstract}

 The lecture explores several problems of analytic number theory in the context of
function fields over a finite field, where they can be approached by
methods different than those of traditional analytic number theory.
The resulting theorems can be used to check existing conjectures
over the integers, and to generate new ones. Among the problems
discussed are:  Counting primes in short intervals and in arithmetic
progressions; Chowla's conjecture on the autocorrelation of the
M\"obius function; and the additive divisor problem.
\end{abstract}




\maketitle

\section{Introduction}

 The goal of this lecture is to explore traditional problems of
analytic number theory in the context of function fields over a
finite field. Several such problems which are currently viewed as
intractable over the integers, have recently been addressed  in the
function field context with vastly different tools than those of
traditional analytic number theory, and the resulting theorems can
be used to check existing conjectures over the integers, and to
generate new ones. The problems that I will address concern
\begin{itemize}
\item Counting primes in short intervals and in arithmetic
progressions
\item
Chowla's conjecture on the autocorrelation of the M\"obius function
\item
The twin prime conjecture
\item
The additive divisor problem
\item The variance of sums of arithmetic
functions in short intervals and arithmetic progressions.
\end{itemize}

Before describing the problems, I will briefly survey some
quantitative aspects of the arithmetic of the ring of polynomials
over a finite field.

\section{{Background on arithmetic in $\fq[x]$  }}

\subsection{The Prime Polynomial Theorem}
Let $\fq$ be a finite field of $q$ elements, and $\fq[x]$ the ring
of polynomials with coefficients in $\fq$. The polynomial  ring
$\fq[x]$ shares several qualitative properties with the ring of
integers $\Z$, for instance having a Euclidean algorithm, hence
unique factorization into irreducibles. There are also several
common quantitative aspects. To set  these up, I review some basics.

The units of the ring of integers are $\pm 1$, and every nonzero
integer is a multiple by a unit of a positive integer. Analogously,
the units of $\fq[x]$ are the nonzero scalars $\fq^\times$, and
every nonzero polynomial is is a multiple by a unit of a monic
polynomial. The analogue of a (positive) prime is a monic
irreducible polynomial. To investigate arithmetic properties of
``typical'' integers, one samples them uniformly in the dyadic
interval $[X,2X]$ with $X\to \infty$; likewise to investigate
arithmetic properties of ``typical'' polynomials, one samples them
uniformly from the monic polynomials $\EM_n$ of degree $n$, with
$\#\EM_n=q^n\to \infty$.

The Prime Number Theorem (PNT) states that the number $\pi(x)$ of
primes $p\leq x$ is asymptotically equal to
\begin{equation}\label{PNT}
\pi(x)\sim \Li(x):=\int_2^x \frac{dt}{\log t}\sim \frac{x}{\log
x}\;,\quad x\to \infty \;.
\end{equation}
The Riemann Hypothesis is equivalent to the assertion that
\begin{equation}\label{RH}
\pi(x)= \Li(x) + O\Big(x^{1/2+o(1)}\Big) \;.
\end{equation}

The Prime Polynomial Theorem asserts that the number $\pi_q(n)$ of
monic irreducible polynomials  of degree $n$ is
\begin{equation}\label{PPT}
\pi_q(n)=\frac{q^n}{n} + O\Big(\frac{q^{n/2}}{n}\Big) \;,
\end{equation}
the implied constant absolute. This corresponds to the PNT (and to
the Riemann Hypothesis) if we map $x\leftrightarrow q^n$, recalling
that $x$ is the number of positive integers up to $x$ and $q^n$ is
the number of monic polynomials of degree $n$. Note that \eqref{PPT}
gives an asymptotic result whenever $q^n\to \infty$; in comparison,
the results described below will usually be valid only in the large
finite field limit, that is $n$ fixed and $q\to \infty$.

\subsection{Cycle structure}

For $f\in \fq[x]$ of positive degree $n$, we say its cycle structure
is $\lambda(f) = (\lambda_1,\dots, \lambda_n)$
if in the prime decomposition $f=\prod_\alpha P_\alpha$ (we allow
repetition), we have $\#\{\alpha : \deg P_\alpha=j\} = \lambda_j$.
In particular  $\deg f = \sum_j j \lambda_j$. Thus we get a
partition of $\deg f$, which we denote by $\lambda(f)$. For
instance, $\lambda_1(f)$ is the number of roots of $f$ in $\fq$, and
 $f$ is totally split in $\fq[x]$ - that is $f(x) = \prod_{j=1}^n
(x- a_j )$, $a_j \in \fq$- if and only if   $\lambda(f) = (n, 0,
\dots,0)$. Moreover $f$ is prime if and only if
$\lambda(f)=(0,0,\dots,0,1)$.

The cycle structure of a permutation $\sigma$ of $n$ letters is
$\lambda(\sigma) = (\lambda_1,\dots, \lambda_n)$ if in the
decomposition of $\sigma$ as a product of disjoint cycles, there are
 $\lambda_j$ cycles of length $j$. For instance, $\lambda_1(\sigma)$
 is the number of fixed points of $\sigma$, and $\sigma=I$ is the identity if and only if $\lambda(\sigma) =
 (n,0,\dots)$. Moreover $\sigma\in S_n$ is an $n$-cycle if and only
 if $\lambda(\sigma)=(0,0,\dots,0,1)$.

For each partition $\lambda \vdash n$, denote by $p(\lambda)$ the
probability that a random permutation on $n$ letters has cycle
structure $\lambda$:
\begin{equation}\label{def of plambda}
p(\lambda) = \frac{\#\{\sigma\in S_n: \lambda(\sigma)=
\lambda\}}{\#S_n } \;.
\end{equation}
  Cauchy's formula for $p(\lambda)$ is
\begin{equation}\label{Cauchy}
p(\lambda) = \prod_{j=1}^n \frac 1{j^{\lambda_j} \cdot \lambda_j!}
\end{equation}
In particular,  the proportion of $n$-cycles in the symmetric group
$S_n$ is $1/n$.

The connection between cycle structures of polynomials and of
permutations is by means of the following observation, a
straight-forward consequence of the Prime Polynomial Theorem
\eqref{PPT}: Given a partition $\lambda\vdash n$, the probability
that a random monic polynomial $f$ of degree $n$ has cycle structure
$\lambda$ is asymptotic, as $q\to \infty$, to the probability
$p(\lambda)$ that a random permutation of $n$ letters has that cycle
structure:
\begin{equation}\label{number of pols with cycle}
\frac 1{q^n} \#\{f \;{\rm monic}, \deg f=n: \lambda(f) = \lambda\} =
p(\lambda) +O\Big( \frac 1q \Big) \;.
\end{equation}
Note that unlike the Prime Polynomial Theorem \eqref{PPT}, this
result \eqref{number of pols with cycle} gives an asymptotic only in
the large finite field limit $q\to \infty$, $n$ fixed.

Having set up the preliminaries, I turn to discussing new results on
quantitative aspects of arithmetic in $\fq[x]$.

\section{Asymptotics in short intervals and arithmetic progressions}

\subsection{Primes in short intervals}
Some of the most important problems in prime number theory concern
the distribution of primes in short intervals and in arithmetic
progressions. According to the Prime Number Theorem, the density of
primes near $x$ is $1/\log x$. Thus one wants to know what is the
number $\pi(x,H)$ of primes in an interval of length $H=H(x)\ll x$
around $x$:
\begin{equation}
\pi(x,H):=\#\{x<p\leq x+H: p\;{\rm prime}\}\;.
\end{equation}
We expect that for $H$ sufficiently large,
\begin{equation}\label{primes short intervals}
\pi(x,H) \sim \frac{H}{\log x} \;.
\end{equation}
The PNT implies that \eqref{primes short intervals} holds for
$H\approx x$, and the Riemann Hypothesis gives \eqref{primes short
intervals} for all $H>x^{1/2+o(1)}$. In 1930, Hoheisel gave an
unconditional proof that \eqref{primes short intervals} holds for
all $H> x^{1-\delta}$ for any positive $\delta<1/33,000$; this has
since been improved, currently to $H>x^{7/12-o(1)}$ (Heath Brown
1988). It is  believed that the result should hold for all
$H>x^\epsilon$, for any $\epsilon>0$, though Maier \cite{Maier}
showed that it does  not  hold for $H=(\log x)^N$ for any $N$; see
Granville and Soundararajan \cite{GS} for a general framework for
such results on irregularities of distribution
 and for sharper results. Selberg (1943) showed, assuming the Riemann Hypothesis, that \eqref{primes short
intervals} holds for {\em almost all} $x$ provided $H/(\log x)^2 \to
\infty$.

To set up an analogous problem for the polynomial ring $\fq[x]$, we
first need to define short intervals.   For a nonzero polynomial
$f\in \fq[x]$, we define its norm by
$$
|f| = \#\fq[x]/(f)= q^{\deg f}\;,
$$
in analogy with the norm of a nonzero integer $0\neq n\in \Z$, which
is $|n|=\#\Z/n\Z$. Given a monic polynomial $A\in \EM_n$ of degree
$n$, and $h<n$, the "short interval" around $A$ of diameter $q^h$ is
the set
\begin{equation}\label{def short interval}
I(A;h):=\{f\in \EM_n: |f-A|\leq q^h\} \;.
\end{equation}
The number of polynomials in this "interval" is
\begin{equation}
H:=\#I(A;h) =q^{h+1} \;.
\end{equation}

 We wish to count the number of prime polynomials in the interval $I(A;h)$.
In the limit $q\to \infty$, Bank, Bary-Soroker and Rosenzweig
\cite{BBSR} give an essentially optimal short interval result:
\begin{thm}
Fix $3\leq h<n$. Then for every monic polynomial $A$ of degree $n$,
the number of prime polynomials $P$ in the interval
$I(A;h)=\{f:|f-A|\leq q^h\}$ about $A$ satisfies
\begin{equation*}
\#\{P\;{\rm prime}, P\in I(A;h)  \} =\frac{H}{n} \left(1+
O_n(q^{-1/2})\right)\;,
\end{equation*}
the implied constant depending only on $n$.
\end{thm}
For irregularities of distribution analogous to Maier's theorem in
the large degree limit $n\to \infty$ ($q$ fixed), see \cite{Thorne}.

For other applications, we will need a version which takes into
account the cycle structure:
\begin{thm}[\cite{BBSR}]\label{thm BBSR}
Fix $n>1$, $3\leq h<n$ and a partition $\lambda\vdash n$. Then 
for any sequence of finite fields $\fq$, and every monic polynomial
$A$ of degree $n$,
$$ \#\{f\in I(A;h):\lambda(f)=\lambda\} = p(\lambda) H\left(1+
O_n(q^{-1/2})\right)\;,
$$
with $p(\lambda)$ as in \eqref{def of plambda}, \eqref{Cauchy}, the
implied constant depending only on $n$.
\end{thm}

\subsection{Primes in  arithmetic progressions}
Dirichlet's theorem states that  any arithmetic progression $n=
A\bmod Q$ contains infinitely many primes provided that $A$ and $Q$
are coprime, and the prime number theorem in arithmetic progressions
states that for fixed modulus $Q$, the number of such primes $p\leq
x$ is
\begin{equation}\label{PNTAP}
\pi(x;Q,A) \sim\frac{\Li(x)}{\phi(Q)} , \quad x\to \infty\;,
\end{equation}
where $\phi(Q)$ is Euler's totient function, the number of residues
coprime to $Q$. The Generalized Riemann Hypothesis (GRH) asserts
that \eqref{PNTAP} continues to hold for moduli as large as
$Q<X^{1/2-o(1)}$. An unconditional version, for almost all
$Q<x^{1/2-o(1)}$, and all $A\bmod Q$,  is given by the
Bombieri-Vinogradov theorem. Going beyond the GRH, the
Elliott-Halberstam conjecture gives a similar statement for $Q$ as
large as $x^{1-\epsilon}$.

For $\fq[x]$, it is a consequence of the Riemann Hypothesis for
curves over a finite field (Weil's theorem) that given a modulus
$Q\in \fq[x]$ of positive degree, and a polynomial $A$ coprime to
$Q$, the number $\pi_q(n;Q,A)$ of primes $P=A\bmod Q$, $P\in \EM_n$
satisfies
$$
\pi_q(n;Q,A) = \frac{\pi_q(n)}{\Phi(Q)} + O(\deg Q \cdot q^{n/2})\;,
$$
where $\Phi(Q)$ is the number of coprime residues modulo $Q$. For
$q\to \infty$, the main term is dominant as long as $\deg Q<n/2$.

Going beyond the Riemann Hypothesis for curves,  Bank, Bary-Soroker
and Rosenzweig \cite{BBSR} show an individual asymptotic continues
to hold for even  larger moduli in the limit $q\to \infty$:
\begin{thm}[\cite{BBSR}]
If $ 1\leq \deg Q\leq n-3$ then
$$\pi_q(n;Q,A)  = \frac{\pi_q(n)}{\Phi(Q)} \left(1+ O_n(q^{-\frac
12})\right) \;.
$$
\end{thm}
This should be considered as an individual version of the
Elliot-Halberstam conjecture. As in the short interval case, they
have a stronger result which takes into account the cycle structure.

\section{Autocorrelations and twisted convolution}\label{sec:conv}

In this section we describe results on the autocorrelation of
various classical arithmetic functions in the function field
context.

\subsection{Autocorrelations of the M\"obius function and Chowla's conjecture}

Equivalent formulations of the PNT and the Riemann Hypothesis can be
given in terms of growth of partial sums of the M\"obius function,
defined by $\mu(n) = (-1)^k$ if $n$ is a product of $k$ distinct
primes, and $\mu(n)=0$ otherwise: The PNT is equivalent to
nontrivial cancellation $\sum_{n\leq x}\mu(n)  = o(x)$, and the RH
is equivalent to square-root cancellation: $\sum_{n\leq x} \mu(n) =
O(x^{1/2+o(1)})$.

A conjecture of Chowla   on the auto-correlation of the M\"obius
function,   asserts that
 given an $r$-tuple of distinct integers $\alpha_1,\dots,\alpha_r$ and
 $\epsilon_i\in \{1,2\}$, not all even,
 then
\begin{equation}\label{chowla's conjecture}
\lim_{N\to \infty} \frac 1N \sum_{n\leq N}\mu(n+\alpha_1
)^{\epsilon_1}\cdot \dots \cdot \mu(n+\alpha_r)^{\epsilon_r}=0 \;.
\end{equation}
Note that the number of nonzero summands here, that is the number of
$n\leq N$ for which $n+\alpha_1,\dots n+\alpha_r$ are all
square-free, is asymptotically $\mathfrak S(\alpha)N$, where
$\mathfrak S(\alpha)>0$ if the numbers $\alpha_1,\dots, \alpha_r$ do
not contain a complete system of residues modulo $p^2$ for every
prime $p$,
so that Chowla's conjecture \eqref{chowla's conjecture} addresses
non-trivial cancellation in the sum.  At this time, the only known
case of Chowla's conjecture \eqref{chowla's conjecture}  is $r=1$
where it is equivalent with the Prime Number Theorem.

 Sarnak \cite{Sarnak} showed that Chowla's conjecture implies that $\mu(n)$  does not correlate
 with  any ``deterministic'' (i.~e., zero entropy) sequence.
For recent studies on the correlation between $\mu(n)$ and several
sequences of arithmetic functions, see \cite{Green Tao, CS, BSZ,
LiuSar}.

In joint work with Dan Carmon \cite{CarmonRudnick}, we have resolved
a version of Chowla's conjecture for $\fq[x]$  in the limit
$q\to\infty$. To formulate it, one defines the M\"obius function of
a nonzero polynomial $F\in \fq[x]$ to be $\mu(F)=(-1)^r$ if
$F=cP_1\dots P_r$ with $0\neq c\in \fq$ and $P_1,\dots, P_r$ are
distinct monic irreducible polynomials, and $\mu(F)=0$ otherwise.

\begin{thm}\label{main thm CR}
Fix $r>1$ and assume that $n>1$ and $q$ is odd.
Then for any choice of distinct polynomials  $\alpha_1,\dots,
\alpha_r \in \fq[x]$, with $\max \deg \alpha_j<n$,  and
$\epsilon_i\in \{1,2\}$, not all even,
\begin{equation}\label{statement of chowla q}
|\sum_{F\in \EM_n} \mu(F+\alpha_1)^{\epsilon_1}\dots
\mu(F+\alpha_r)^{\epsilon_r}|  \ll_{r,n} q^{n-\frac 12} \;.
\end{equation}
\end{thm}
Thus for fixed $r,n>1$,
\begin{equation}
\lim_{q\to \infty} \frac 1{\#\EM_n}\sum_{F\in \EM_n}
\mu(F+\alpha_1)^{\epsilon_1}\dots \mu(F+\alpha_r)^{\epsilon_r} =0
\end{equation}
under the assumption of Theorem~\ref{main thm CR}, giving an
analogue of Chowla's conjecture \eqref{chowla's conjecture}.

Note that the number of square-free  monic polynomials of degree $n$
is, for $n>1$, equal to $q^n-q^{n-1}$.
 Hence, given $r$  distinct polynomials $\alpha_1,\dots,
\alpha_r \in \fq[x]$, with $\deg \alpha_j<n$, the number of $F\in
\EM_n$ for which all of $F(x)+\alpha_j(x)$ are square-free is $q^n +
O(rq^{n-1}) $ as $q\to \infty$. Thus indeed we display cancellation.

The starting point in our argument is Pellet's formula,
which asserts that for the polynomial ring $\fq[x]$ with $q$ odd,
the M\"obius function $\mu(F)$ can be computed in terms of the
discriminant $\disc(F)$ of $F(x)$ as
\begin{equation}\label{Pellet}
\mu(F) = (-1)^{\deg F} \chi_2(\disc(F))\;,
\end{equation}
where $\chi_2$ is the quadratic character of $\fq$. That   allows us
to express the LHS of \eqref{statement of chowla q}  as an
$n$-variable character sum and to estimate it by freezing all but
one of the variables, and then using the Riemann Hypothesis for
curves (Weil's theorem) to bound the one-variable sum. A key point
is to bound the number of times when there is no cancellation in the
one-variable sum.


\subsection{Twin primes}
 It is an ancient conjecture that there are infinitely many twin
primes, and a refined quantitative form, due to Hardy and
Littlewood, asserts that given distinct integers $a_1$, $\dots$,
$a_r$, the number $\pi(x; a_1,\dots, a_r)$ of integers $n\leq x$ for
which $n+a_1,\dots ,n+a_r$ are simultaneously prime is
asymptotically
\begin{equation}\label{HL conj}
\pi(x;a_1,\dots,a_r) \sim \mathfrak S(a_1,\dots,a_r)\frac x{(\log
x)^r}, \quad x\to \infty\;,
\end{equation}
for a certain constant $\mathfrak S(a_1,\dots, a_r)$, which is
positive whenever there are no local congruence obstructions.
Despite the striking recent breakthroughs by Zhang \cite{Zhang} and
Maynard \cite{Maynard}, this conjecture is still open even for $r=2$
(twin primes).

Recently the function field version of the problem was solved.
 Bary-Soroker \cite{BSIMRN} proved that for given $n,r$
then for any sequence of finite fields $\fq$ of odd cardinality $q$,
and distinct polynomials $a_1,\dots,a_r\in \fq[x]$ of degree less
than $n$, the number $\pi_q(n;a_1,\dots,a_r)$ of monic polynomials
$f\in \fq[x]$ of degree $n$ such that $f+a_1,\dots, f+a_r$ are
simultaneously irreducible satisfies
\begin{equation}
\pi_q(n;  a_1,\dots, a_r)\sim \frac{q^n}{n^r}, \quad q\to \infty \;.
\end{equation}
This improves on earlier results by Pollack \cite{Pollack} and by
Bary-Soroker \cite{BSadv}.

\subsection{The additive divisor problem}
The divisor function $\divid_r(n)$ is the number of ways of writing
a positive integer $n$ as a product of $r$ positive integers. In
particular for $r=2$ we recover the classical divisor function
$\divid_2(n) = \sum_{d\mid n}1$. The mean value of $\divid_r$ is
\begin{equation}
\frac 1x \sum_{n\leq x} \divid_r(n) \sim \frac {(\log
x)^{r-1}}{(r-1)!},\quad x\to \infty\;.
\end{equation}
Likewise, the divisor function $\divid_r(f)$  for a  monic
polynomial $f\in \fq[x]$ is defined as the number of $r$-tuples of
monic polynomials $(a_1,\dots,a_r)$ so that $f=a_1\cdot \ldots \cdot
a_r$. The mean value of $\divid_r$, when averaged over all monic
polynomials of degree $n$,  is
\begin{equation}
 \frac 1{q^n} \sum_{f\in \EM_n} \divid_r(f) = \binom{n+r-1}{r-1}  =
 \frac{n^{r-1}}{(r-1)!}+\dots\;,
\end{equation}
which is a polynomial of degree $r-1$ in   $n$.


The "additive divisor problem" (other names are "shifted divisor"
and "shifted convolution") is to understand the autocorrelation of
the divisor function, that is  the sum (where $h\neq 0$ is fixed for
this discussion)
\begin{equation}
D_r(X;h):=\sum_{n\leq X}\divid_r(n)\divid_r(n+h)\;.
\end{equation}
These sums are of importance in studying the moments of the Riemann
$\zeta$-function on the critical line, see \cite{Ivic, CG}.

For $r=2$ (the ordinary divisor function), Ingham \cite{Ingham}  and
Estermann \cite{Estermann} showed that
\begin{equation}
\sum_{n\leq X} \divid_2(n)\divid_2(n+h) \sim X P_2(\log X;h),\quad
X\to \infty
\end{equation}
where $P_2(u;h)$ is a quadratic polynomial in $u$.



For $r\geq 3$ it is conjectured that
\begin{equation}
D_r(X;h)\sim XP_{2(r-1)}(\log X;h),\quad X\to \infty
\end{equation}
where  $P_{2(r-1)}(u;h)$ is a polynomial in $u$ of degree $2(r-1)$,
whose coefficients depend on $h$ (and $r$). However to date one is
very far from being able to even get good upper bounds  on
$D_r(X;h)$.
Moreover, even a conjectural description of the polynomials
$P_{2(r-1)}(u;h)$ is difficult to obtain, see \cite{Ivic, CG}.


In joint work with Andrade and Bary-Soroker \cite{ABSR},  we study a
version of the additive divisor problem for $\fq[x]$.
 We show:
\begin{thm}\label{main thm}
Let $0\neq h\in \fq[x] $, and $n>\deg h$. Then for $q$ odd,
\begin{equation}
\frac 1{q^n} \sum_{f\in \EM_n} \divid_r(f)\divid_r (f+h) =
\binom{n+r-1}{r-1}^2  + O_n(  q^{-1/2})\;,
\end{equation}
the implied constant depending only on $n$.
\end{thm}
Note that $\binom{n+r-1}{r-1}^2$ is a polynomial in $n$ of degree
$2(r-1)$ with leading coefficient $1/[(r-1)!]^2$.



\subsection{About proofs}
The results of this section can all be deduced from one principle
(though this was not the original proof of most), namely that for a
random monic polynomial $f\in \EM_n$ of degree $n$, the cycle
structure of $f$ and its shift $f+\alpha$ are {\em independent} as
$q\to \infty$. Precisely, in \cite{ABSR} we show that for for fixed
$n>1$, and two partitions $\lambda',\lambda''\vdash n$, given any
sequence of finite fields $\fq$ of odd  cardinality $q$, and nonzero
$\alpha\in \fq[x]$ of degree less than $n$, then
\begin{equation}\label{independence}
\lim_{q \to \infty} \frac 1{q^n} \#\{f\in \EM_n:\lambda(f) =
\lambda',\lambda(f+\alpha) = \lambda''\} = p(\lambda')\times p(
\lambda'')
\end{equation}
where $p(\lambda)$, as in \eqref{def of plambda}, \eqref{Cauchy},
is the probability that a random permutation on $n$ letters has
cycle structure $\lambda$.
This result  
is an elaboration of earlier work by Bary-Soroker \cite{BSIMRN}
  which dealt with the case of
$n$-cycles, where $\lambda=\tilde \lambda=(0,\dots,0,n)$. There is
also a version allowing several distinct shifts.

To prove  \eqref{independence} we need to compute a certain Galois
group: Let $\FF$ be an algebraic closure of $\fq$,   $\mathbf A =
(A_0,\dots,A_{n-1})$ be indeterminates, and
\begin{equation}
\mathcal F (\mathbf A,x) = x^n+A_{n-1}x^{n-1} + \dots +A_0
\end{equation}
the generic polynomial of degree $n$, whose Galois group over
$\FF(\mathbf A)$ is well-known to be the full symmetric group $S_n$.
For nonzero $\alpha\in \fq[x]$ of degree less than $n$, let
\begin{equation}
\mathcal G(\mathbf A,x) =\mathcal F (\mathbf A,x)\Big(\mathcal
F(\mathbf A,x) + \alpha(x)\Big) \;.
\end{equation}
Bary-Soroker \cite{BSIMRN} shows that for odd $q$, the Galois group
of $\mathcal G$ over $\FF(\mathbf A)$ is  the product $S_n\times
S_n$, the maximal possible group. The proof requires an ingredient
from the proof of Chowla's conjecture \cite{CarmonRudnick} discussed
above.

Once we know the Galois group of $\mathcal G(\mathbf A,x)$, we apply
an explicit version of Chebotarev's theorem for function fields to
prove \eqref{independence}, see \cite{ABSR} for the details.

\section{The variance of sums of arithmetic functions 
and matrix integrals} \label{sec:variance}

I now describe some  results concerning the variance of sums of
several arithmetic functions. A common feature is that the variance
is expressed as a matrix integral.

\subsection{Variance of primes in short intervals}

The von Mangoldt function is defined as $\Lambda(n)=\log p$ if
$n=p^k$ is a prime power, and $0$ otherwise. A form of the Prime
Number Theorem (PNT) is the assertion that
\begin{equation}
\psi(x):=\sum_{n\leq x}\Lambda(n) \sim x \quad \mbox{ as }x\to
\infty \;.
\end{equation}

To study  the distribution of primes in short intervals, we define
for  $1\leq H\leq x$,
\begin{equation}
 \psi(x;H):=   \sum_{n\in[x-\frac H2,x+\frac H2]} \Lambda(n) \;.
\end{equation}
The Riemann Hypothesis guarantees an asymptotic formula
$\psi(X;H)\sim H$ as long as $H>X^{\frac 12 +o(1)}$.
Goldston and Montgomery \cite{GM} studied the variance of
$\psi(x;H)$, relating it to the pair correlation function of the
zeros of the Riemann zeta function.
The conjecture of Goldston and Montgomery, as refined by Montgomery
and Soundararajan\footnote{based on Hardy-Littlewood type heuristics
} \cite{MSrefined}  is that in the range $X^\epsilon
<H<X^{1-\epsilon}$, as $X\to \infty$:
\begin{equation} \label{GMrefined}
\frac 1X \int_1^X \left |\psi(x;H)-H\right|^2 dx \sim H\Big(\log
X-\log H-(\gamma + \log 2\pi)\Big)
\end{equation}
 with $\gamma$ being Euler's constant.

With J.~Keating, we prove a function field analogue of
Conjecture~\ref{GMrefined}:
 \begin{thm}[\cite{KR}]
 For $h\leq n-5$, as $q\to \infty$,
\begin{equation*}
\frac 1{q^n} \sum_{A\in \EM_n} \Big|\sum_{|f-A|\leq q^h} \Lambda(f)
-H \Big|^2 \sim H\int_{U(n-h-2)} \Big| \tr U^n \Big|^2dU=H
(n-h-2)\;.
\end{equation*}
\end{thm}
Recall $H:=\#\{f:|f-A|\leq q^h\} = q^{h+1}$. Here the matrix
integral is over the  unitary group $U(n-h-2)$, equipped with its
Haar probability measure.

\subsection{Variance of primes in arithmetic progressions}
  A form of the Prime Number Theorem for arithmetic progression states that for a modulus
$Q$ and $A$ coprime to $Q$,
\begin{equation}\label{PNT for arith prog}
   \psi(X;Q,A) := \sum_{\substack{n\leq X\\ n=A\bmod Q}} \Lambda(n)
   \sim \frac{X}{\phi(Q)},\quad \mbox{ as }X\to \infty
   \;.
\end{equation}

In most arithmetic applications it is crucial to allow the modulus
to grow with $X$.
 For very large moduli $Q>X$, there can be at
most one prime in the arithmetic progression $P=A \mod Q$ so that
the interesting range is $Q<X$. To study the fluctuations of
$\psi(X;Q,A)$,  define
\begin{equation}
  G(X,Q)=\sum_{\substack{A\bmod Q\\ \gcd(A,Q)=1}} \left| \psi(X;Q,A)-\frac X{\phi(Q)}
  \right|^2 \;.
\end{equation}
Hooley, in his ICM article  \cite{HooleyICM},  conjectured that
under some (unspecified) conditions,
\begin{equation}\label{Hooley conj}
  G(X,Q) \sim X\log Q \;.
\end{equation}
Friedlander and Goldston \cite{FG} conjecture that \eqref{Hooley
conj} holds if $X^{1/2+\epsilon}<Q<X$, and further conjecture that
if $X^{1/2+\epsilon}<Q<X^{1-\epsilon}$ then
\begin{equation}\label{FG conj}
  G(X,Q)= X \left( \log Q -\Big(\gamma +\log 2\pi +
\sum_{p\mid Q} \frac{\log p}{p-1} \Big) \right)+o(X) \;.
\end{equation}
They show  that both \eqref{Hooley conj} (in the range
$X^{1/2+\epsilon}<Q<X$) and \eqref{FG conj} (in the range
$X^{1/2+\epsilon}<Q<X^{1-\epsilon}$) hold assuming GRH and a strong
version of the Hardy-Littlewood conjecture \eqref{HL conj} on prime
pairs. For $Q<X^{1/2}$   little is known. In any case, Hooley's
conjecture \eqref{Hooley conj} has not been proved in any range.

With J.~ Keating \cite{KR} we resolve the function-field version  of
Conjecture~\eqref{Hooley conj}:

\begin{thm}\label{thm KR primes in AP}

Fix $n\geq 2$. Given a sequence of finite fields $\fq$ and
square-free polynomials $Q(x)\in \fq[x]$  with $2\leq \deg Q\leq
n-1$,  then as  $ q\to \infty$,
\begin{equation}\label{RMT conseq}
    G(n;Q) \sim q^n\int_{U(\deg Q-1)} |\tr U|^ndU =q^n(\deg Q-1) \;.
\end{equation}
 \end{thm}


We can compare our result \eqref{RMT conseq} to the conjectures
 \eqref{Hooley conj} and \eqref{FG conj}:
The range $X^{1/2}<Q<X^{}$ corresponds to $\deg Q<n<2\deg Q$, so
that we recover the function field version of conjecture
\eqref{Hooley conj}; note that \eqref{RMT conseq} holds for all $n$,
not just in that range. Thus we believe that Hooley's conjecture
\eqref{Hooley conj} should hold for all $Q>X^{\epsilon}$. We refer
to Fiorilli's recent work \cite{Fiorilli} for a more refined
conjecture in this direction.

 \subsection{Almost-primes}
A variation on this theme was proposed by B.~Rodgers \cite{Rodgers}.
Instead of primes, he considered "almost primes", that is products
of two prime powers. A useful weight function for these is the
generalized von Mangoldt function
\begin{equation}
\Lambda_2 =\Lambda\conv \Lambda + \deg \cdot \Lambda = \mu\conv
\deg^2
\end{equation}
which is supported on products of two prime powers ($\conv$ means
Dirichlet convolution). The mean value of $\Lambda_2$ over  the set
$\mathcal M_n$ of monic polynomials of degree $n$ is
\begin{equation}
\frac 1{q^n} \sum_{f\in \mathcal M_n} \Lambda_2(f) =
n^2-(n-1)^2=2n-1 \;.
\end{equation}

To count almost primes in the short intervals, set for
 $A\in \mathcal M_n$, and $1\leq h< n$
 \begin{equation}
  \Psi_2(A;h) = \sum_{ f\in I(A;h) } \Lambda_2(f)\;.
\end{equation}
 Rodgers showed \cite{Rodgers} that the variance of $\Psi_2(A;h)$ is given as $q\to \infty$,
for fixed $n$ and $h\leq n-5$, by the matrix integral
\begin{equation}
  \Var\Psi_2(\bullet;h) \sim H  \int_{U(n-h-2)} \Big|
\sum_{j=1}^{n-1} \tr U^j\tr U^{n-j}  -n\tr U^n \Big|^2 dU,\quad q\to
\infty\;.
\end{equation}
He shows the matrix integral to be equal to
$(4(n-h-2)^3-(n-h-2))/3$, in fact that
\begin{equation}
 \int_{U(N)} \Big|
\sum_{j=1}^{n-1} \tr U^j\tr U^{n-j}  -n\tr U^n \Big|^2 dU
=\sum_{d=1}^{\min(n,N)}(d^2-(d-1)^2))^2 \;.
\end{equation}

\subsection{Sums of the M\"obius function and the Good-Churchhouse conjecture}

It is a standard heuristic to assume that the M\"obius function
behaves like a  random variable taking values $\pm 1$ with equal
probability, and supported on the square-free integers  (which have
density $1/\zeta(2) = 6/\pi^2$). In particular if we consider the
sums of $\mu(n)$ in blocks of length $H$,
\begin{equation}
M(x;H):=\sum_{|n-x|<H/2} \mu(n)
\end{equation}
then   when averaged over $x$,
$M(x,H)/\sqrt{H}$ has mean zero, 
and it was conjectured by Good and Churchhouse \cite{GC} in 1968
that $M(x;H)/\sqrt{H}$ has variance $1/\zeta(2)$:
\begin{equation}\label{GC conj}
\frac 1X\int_X^{2X} |M(x;H)|^2 \sim \frac H{\zeta(2)}
\end{equation}
for $X^\epsilon<H=H(X)<X^{1-\epsilon}$. Moreover they conjectured
that the normalized sums $M(x;H)/\sqrt{H/\zeta(2)}$ have
asymptotically a normal distribution.

We can apply our method to evaluate the variance of sums of the
M\"obius function in short intervals for $\fq[x]$.   Set
\begin{equation}
\EN_\mu(A;h):= \sum_ {f\in I(A;h)} \mu(f) \;.
\end{equation}
The mean value of $\EN_\mu(A;h)$   is  $0$, and the variance is
\begin{thm}[Keating-Rudnick \cite{KR14}]\label{thm:Mobius}
If  $h\leq n-5$ then as $q\to \infty$,
\begin{equation*}
\Var\EN_\mu(\bullet;h)\sim H\int_{U(n-h-2)}| \tr \Sym^n U|^2 dU = H
\end{equation*}
where $\Sym^n$ is the representation of the unitary group $U(N)$ on
polynomials of degree $n$ in $N$ variables.
\end{thm}
Theorem~\ref{thm:Mobius}  is consistent with Conjecture~\eqref{GC
conj} if we replace $H$ by  $H/\zeta_q(2)$ where $\zeta_q(2) =
\sum_{f}  1/|f|^2$ (the sum over all monic $f$), which tends to $1$
as $q\to \infty$.

%

 \subsection{The divisor function in short intervals}

Dirichlet's divisor problem addresses the size of the remainder term
$\Delta_2(x)$ in partial sums of the divisor function:
\begin{equation}\label{def of Delta2}
\Delta_2(x):=\sum_{n\leq x} \divid_2(n)-x\Big(\log
x+(2\gamma-1)\Big)
\end{equation}
 where $\gamma$ is the Euler-Mascheroni constant. For the higher divisor
 functions one defines a remainder term $\Delta_k(x)$ similarly as the
difference between the partial sums $\sum_{n\leq x} \divid_k(n)$ and
a smooth term $xP_{k-1}(\log x)$ where $P_{k-1}(u)$ is a certain
polynomial of degree $k-1$.

Let
\begin{equation}\label{Def of Delta(x:H)}
\Delta_k(x;H) = \Delta_k(x+H)-\Delta_k(x)
\end{equation}
be the remainder term for sums of $d_k$ over short intervals
$[x,x+H]$. Jutila \cite{Jutila1984}, Coppola and Salerno \cite{CS},
and Ivi\'c \cite{Ivic Bordeaux 2009,  Ivic RramanujanJ} show that,
for $X^\epsilon<H< X^{1/2-\epsilon}$,  the mean square of
$\Delta_2(x,H)$ is asymptotically equal to
\begin{equation}\label{Ivic}
\frac 1X\int_X^{2X} \Big(\Delta_2(x,H)\Big)^2 dx \sim HP_3(\log
X-2\log H)
\end{equation}
for a certain cubic polynomial $P_3$.

Lester and Yesha \cite{LY} showed that $\Delta_2(x,H)$, normalized
to have unit mean-square using \eqref{Ivic}, has a Gaussian value
distribution at least for a narrow range of $H$ below $X^{1/2}$:
$H=\sqrt{X}/L$, where $L=L(X)\to\infty $ with $X$, but  $L \ll
X^{o(1)}$,
(see \cite{LY} for the precise statement), the conjecture being that
this should hold for   $X^\epsilon<H<X^{1/2-\epsilon}$ for any
$\epsilon>0$.


In joint work with J. Keating, B. Rodgers and E. Roditty-Gershon
\cite{KRR}, we study the corresponding problem of the sum of
$\divid_k(f)$ over short intervals  for $\fq[x]$. Set
\begin{equation}
\EN_{\divid_k}(A;h):=\sum_{f\in I(A;h)} \divid_k(f) \;.
\end{equation}
The mean value is 
\begin{equation}
\frac 1{q^n}\sum_{A\in \EM_n} \EN_{\divid_k} (A;h)
=q^{h+1}\binom{n+k-1}{k-1} \;.
\end{equation}

In analogy with \eqref{def of Delta2}, \eqref{Def of Delta(x:H)} we
set
\begin{equation}
\Delta_k(A;h):=\EN_{\divid_k}(A;h)-q^{h+1}\binom{n+k-1}{k-1} \;.
\end{equation}
It can be shown that $\Delta_k(A;h)\equiv 0$ vanishes identically
for $h>(1-\frac 1k)n-1$. Using Theorem~\ref{thm BBSR}  \cite{BBSR},
we can show that for all $3\leq h<n$
\begin{equation}
\Delta_k(A;h) \ll_{n,k} q^{h+\frac 12}
\end{equation}
is smaller than the main term.

We express  the mean square of $\Delta_k(A,h)$ (which is the
variance  of $\EN_{\divid_k}(A;h)$) in terms of a matrix integral.
  Let $\Lambda^j:U(N)\to GL(\Lambda^j \C^N)$ be the exterior $j$-th
power representation ($0\leq j\leq N$).  Define the matrix integrals
over the group $U(N)$ of $N\times N$ unitary matrices
\begin{equation}\label{def of I}
I_k(m;N):=\int_{U(N)} \Big| \sum_{\substack{j_1+\dots+ j_k=m\\
0\leq j_1,\dots, j_k \leq N}} \tr\Lambda^{j_1}(U) \dots
\tr\Lambda^{j_k}(U) \Big|^2 dU\;,
\end{equation}
the integral with respect to the Haar probability measure.

By definition, $I_k(m;N)=0$ for $m>kN$. We have a functional
equation $I_k(m;N) = I_k(kN-m;N)$ and
\begin{equation}\label{eq:DG}
I_k(m;N) = \binom{m+k^2-1}{k^2-1},\quad m\leq N\;.
\end{equation}
The identity \eqref{eq:DG} can be proved by various means, for
instance using the work of Diaconis and Gamburd \cite{DG} relating
matrix integrals to counting magic squares.

\begin{thm}[\cite{KRR}]
Let $n\geq 5$, and $h\leq \min(n-5,(1-\frac 1k)n-2)$. Then as $q\to
\infty$,
\begin{equation*}
\frac 1{q^n}\sum_{A\in \EM_n} |\Delta_k(A;h)|^2 \sim  H\cdot
I_k(n;n-h-2)\;.
 \end{equation*}
\end{thm}

In particular  for the standard divisor function ($k=2$), 
if $h\leq n/2-2$ and $n\geq 8$ then
\begin{equation}
\frac 1{q^n}\sum_{A\in \EM_n} |\Delta_2(A;h)|^2 \sim  H \frac{
(n-2h+5) (n-2h+6) ( n-2h+7)}6 \;.
\end{equation}
 This is consistent with  \eqref{Ivic}, which  leads us to expect
a cubic polynomial in $(  n-2h)$.


\section{How to compute the variance}

Our results on variance described in \S~\ref{sec:variance} depend on
expressing the variance in terms of zeros of Dirichlet L-functions
for $\fq[x]$, and using recent equidistribution results of Katz
\cite{KatzKR1}, \cite{KatzKR2}, tailor-made for this purpose. To
describe how this is done, we give some background on L-functions.

\subsection{Dirichlet L-functions}

Let $Q(x)\in \fq[x]$ be a polynomial of positive degree. A Dirichlet
character modulo $Q$ is a homomorphism $\chi:\left( \fq[x]/(Q)
\right)^\times \to \C^\times$. A Dirichlet character  $\chi$ is
``even'' if $\chi(cF)=\chi(F)$ for all $0\neq c\in \fq$, and $\chi$
is {\em primitive} if there is no proper divisor $Q'\mid Q$ so that
$\chi(F)=1$ whenever $F$ is coprime to $Q$ and $F=1\mod Q'$. The
number of Dirichlet characters modulo $Q$ is $\Phi(Q)$, and the
number of even characters modulo $Q$ is
$\Phi^{ev}(Q)=\Phi(Q)/(q-1)$.

The L-function $\mathcal L(u,\chi)$ attached to $\chi$ is defined as
\begin{equation}\label{Def of L}
\mathcal L(u,\chi) = \sum_{\substack{f\;{\rm monic}\\
(f,Q)=1}} \chi(f)u^{\deg f} = \prod_{P\nmid Q} (1-\chi(P)u^{\deg
P})^{-1}
\end{equation}
where the product, over all monic irreducible polynomials in
$\fq[x]$, is absolutely convergent for $|u|<1/q$.


If $Q\in \fq[x]$ is a polynomial of degree $\deg Q\geq 2$, and
$\chi\neq \chi_0$ is a nontrivial character mod $Q$, then the
L-function $\mathcal L(u,\chi)$ is  a polynomial in $u$ of degree at
most $\deg Q-1$.
 Moreover, if $\chi$ is an  even  character
there is a "trivial" zero at $u=1$.


For a primitive even character modulo $Q$, we can write
\begin{equation}\label{frobenius}
\mathcal L(u,\chi) = (1-u) \det (I-uq^{1/2}\Theta_\chi)
\end{equation}
where the matrix $\Theta_\chi\in U(\deg Q-2)$ is unitary (as follows
from the Riemann Hypothesis for curves), uniquely defined up to
conjugacy. It is called the unitarized Frobenius matrix of $\chi$.
Likewise, if $\chi$ is odd and primitive then $\mathcal L(u,\chi) =
\det (I-uq^{1/2}\Theta_\chi)$ where $\Theta_\chi\in U(\deg Q-1)$ is
unitary.

Katz \cite{KatzKR2} showed that as $\chi$ varies over all primitive
even characters modulo $x^{N+2}$, the unitarized Frobenii
$\Theta_\chi$ become uniformly distributed in the projectivized
unitary group $PU(N)$ for $N\geq 3$ as $q\to \infty$ (and also for
$N=2$ if $q$ is coprime to $2$ and $5$). Thus for any nice class
function $F$ on $U(N)$, which is invariant under the center
($F(zU)=F(u)$, $z$ on the unit circle), we have
\begin{equation}\label{Katz even}
\lim_{q \to \infty} \frac 1{\Phi_{ev}(x^{N+2})} \sum_{\substack{\chi
\bmod x^{N+2}\\ {\rm even \; primitive} }} F(\Theta_\chi)  =
\int_{PU(N)}F(U)dU \;.
\end{equation}

\subsection{Short intervals as arithmetic progressions}
Our method to handle sums over short intervals $I(A;h)=\{f:|f-A|\leq
q^h\} $ is to relate them to arithmetic progressions modulo
$x^{n-h}$.

Denote by $\mathcal P_{\leq n}$ the set of all polynomials of degree
at most $n$. We define a map $\inv_n: \mathcal P_{\leq n}\to
\mathcal P_{\leq n}$ by
\begin{equation}
\theta_n(f) = x^n f(\frac 1x)
\end{equation}
 which takes $f(x) = f_0+f_1x+\dots
+f_n x^n$, $n=\deg f$ to the ``reversed'' polynomial
\begin{equation}
  \inv_n(f)(x) = f_0x^n+f_1x^{n-1}+\dots +f_n\;.
\end{equation}
Then for  $B\in \mathcal M_{n-h-1}$, the  map $\inv_n$ takes the
"interval" $I(T^{h+1}B;h)$ bijectively onto the arithmetic
progression  $\{ g\in\mathcal P_{\leq n}: g\equiv \theta_{n-h-1}(B)
\mod x^{n-h} \} $.

\subsection{A formula for the variance}
The identification of short intervals with arithmetic progressions
allows us to express sums of several arithmetic functions in terms
of even Dirichlet characters. For the case of the von Mangoldt
function, this is done in \cite{KR}. I illustrate this
identification in the case of the M\"obius function
(Theorem~\ref{thm:Mobius}): We denoted by $\EN_\mu(A;h) = \sum_{f\in
I(A;h)} \mu(f)$. Then for $B\in \EM_{m-h-1}$,
\begin{equation}\label{centered gamma}
\EN_\mu(T^{h+1}B;h)=
 \frac{1}{\Phi_{ev}(x^{n-h})} \sum_{\substack{\chi \bmod
x^{n-h}\\\chi\neq \chi_0\mbox{ even}}} \bar\chi(\inv_{n-h-1}(B))
(\EM(n;\mu\chi)  -\EM(n-1;\mu\chi))
\end{equation}
where
\begin{equation}
\EM(n;\mu\chi)=\sum_{f\in \mathcal M_n}\mu(f)\chi(f) \;.
\end{equation}

We  next express the sums $\EM(n;\mu\chi)$ in terms of zeros of the
L-function $\mathcal L(u,\chi)$; for $\chi$  primitive this means in
terms of the unitarized Frobenius matrix $\Theta_\chi$. The
connection is made by writing the generating function identity
\begin{equation}
\sum_{n=0}^\infty \EM(n;\mu\chi)u^n = \frac 1{\mathcal L(u,\chi)}\;.
\end{equation}
Therefore we find that for $\chi$ primitive and even,
\begin{equation}\label{formula for M_n}
\EM(n;\mu\chi) = \sum_{k=0}^n q^{k/2} \tr \Sym^k \Theta_\chi
\end{equation}
where for $N>1$, $\Sym^n: GL(N,\C)\to \Sym^n \C^N$ is the symmetric
$n$-th power representation.
Consequently we obtain
\begin{equation}
\Var  \EN_\mu(\bullet;h) =  \frac {q^{h+1}}   {\Phi_{\rm
ev}(x^{n-h})}
\sum_{\substack{\chi \bmod x^{n-h}\\
\chi {\rm \; even \; and\;  primitive}}} |\tr \Sym^n \Theta_\chi |^2
+O(q^h) \;.
\end{equation}

Using Katz's equidistribution theorem \eqref{Katz even} we get
\begin{equation}
\lim_{q\to \infty}\frac{\Var(\EN_\mu(\bullet;h))} {q^{h+1}} =
\int_{PU(n-h-2)} \left|\tr \Sym^n U \right|^2 dU \;.
\end{equation}
The matrix integrals equals $1$, hence
we  conclude that $\Var(\EN_\mu(;h))\sim q^{h+1}=H$, which is
Theorem~\ref{thm:Mobius}.

\section{Acknowledgements} I am grateful to J.~Andrade, L.~Bary-Soroker, J.~Keating, E. Kowalski, S.~Lester,
  E.~Roditty-Gershon and K.~Soundararajan for their comments on
earlier versions of this survey.





\end{document}